\documentclass[12pt,reqno]{amsart}
\usepackage{amsmath,amssymb}

\newcommand{\svskip}{\vspace{3mm}}

\newcommand{\C}{{\mathbb C}}
\newcommand{\Q}{{\mathbb Q}}
\newcommand{\Z}{{\mathbb Z}}
\newcommand{\A}{{\mathbb A}}

\newcommand{\fm}{{\mathfrak m}}

\newcommand{\Proof}{\noindent{\bf Proof.}\quad}
\newcommand{\Spec}{{\rm Spec}\:}
\newcommand{\Ker}{{\rm Ker}\:}

\newcommand{\transdeg}{{\rm trans.deg}}

\newcommand{\GL}{{\rm GL}}

\newcommand{\wt}{\widetilde}
\newcommand{\ol}{\overline}

\newcommand{\QED}{{\unskip\nobreak\hfil\penalty50\quad\null\nobreak\hfil
{$\Box$}\parfillskip0pt\finalhyphendemerits0\par\medskip}}

\newtheorem{thm}{Theorem}[section]
\newtheorem{lem}[thm]{Lemma}
\newtheorem{prop}[thm]{Proposition}

\newtheorem{example}[thm]{Example}

\begin{document}
\title{Algebraic derivations on affine domains}
\dedicatory{To Professor R.V. Gurjar on his sixtieth birthday}
\author{Kayo Masuda \and Masayoshi Miyanishi}
\thanks{Supported by Grant-in-Aid for Scientific Research (C) 
22540059 \& 21540055, JSPS} 
\address{School of Science and Technology \\
Kwansei Gakuin University \\
2-1 Gakuen, Sanda, 669-1337, Japan}
\email{kayo@kwansei.ac.jp}

\smallskip

\address{Research Center for Mathematical Sciences \\
Kwansei Gakuin University \\
2-1 Gakuen, Sanda, 669-1337, Japan}
\email{miyanisi@kwansei.ac.jp}

\date{}

\maketitle

\begin{abstract}
We define an {\em algebraic} derivation on an affine domain $B$ 
defined over an algebraically closed  field $k$ of characteristic $0$, 
which is called a {\em locally finite} derivation in \cite{Essen}, for example, and has appeared in commutative and 
non-commutative contexts in other references. 
We, without being aware of this existing definition and related results, 
introduced the term of algebraic derivation by extracting a property analogous to algebraic actions of 
algebraic groups. The first section is devoted to the graded ring structure which the algebraic derivation 
$D$ defines on $B$ in a natural fashion. The graded ring structure is indexed by an abelian monoid which is 
a submonoid of the additive group of the ground field $k$. This structure is already observed in \cite{K}. 
But our approach is more computational and straightforward. If the monoid indexing the graded ring structure of $B$ is 
rather restricted (see Theorem \ref{Theorem 1.9} below), the derivation $D$ is close to what is called an {\em Euler derivation} mixed with a locally nilpotent derivation. We observe this fact when $B$ is a polynomial ring mostly in dimension two. In fact, the results in  section two give various characterization of a polynomial ring $k[x,y]$ in terms of algebraic derivations. 
The third section gives a remark on singularities which can coexist with algebraic derivations. The results given in sections two and three are new.
\end{abstract}

\section*{Introduction} 

Let $G$ be an algebraic group defined over an algebraically closed field $k$. If $G$ acts on an affine algebraic 
variety $X=\Spec B$, we say that the action is {\em algebraic} if for every element $b$ of $B$ the $k$-vector 
subspace $\sum_{g\in G(k)}~k\cdot{}^gb$ in the ring $B$ generated by all translates ${}^gb$ with $g \in G(k)$ 
has finite dimension, where $G(k)$ is the set of closed points of $G$. Suppose that $G$ is an affine algebraic 
group $\Spec R$ and the $G$-action is given by a $k$-morphism $\sigma : G \times X \to X$. Let 
$\sigma^\ast : B \to R\otimes B$ be the coaction. For $b \in B$, write 
\[
\sigma^\ast(b)=\sum_{i=1}^n f_i\otimes b_i, \quad f_i \in R,~~ b_i\in B.
\]
Then ${}^gb=\sum_{i=1}^n\pi_g(f_i)\cdot b_i$, where $\pi_g : R \to k$ is the $k$-algebra homomorphism corresponding to 
a point $g \in G(k)$. Then $\sum_{g \in G(k)}~k\cdot{}^gb$ is a $k$-subspace of $\sum_{i=1}^nk\cdot b_i$, whence
$\sum_{g\in G(k)}~k\cdot{}^gb$ has finite dimension. 

Suppose that the ground field $k$ has characteristic zero and 
$G$ is the additive group scheme $G_a$. 
Then a $G_a$-action on an affine algebraic variety $X=\Spec B$ corresponds 
bijectively to a locally nilpotent $k$-derivation $\delta$ on $B$. 
In this case the $k$-subspace $\sum_{g \in G_a(k)}~k\cdot{}^gb$ coincides with 
$\sum_{n\ge 0}k\cdot\delta^n(b)$, which has finite dimension because $\delta^n(b)=0$ for $n \gg 0$. Note that if $\xi$ 
is an element of the quotient field $Q(B)$, the $k$-vector space $\sum_{n \ge 0}k\cdot\delta^n(\xi)$ does not necessarily 
have finite dimension. An example is a locally nilpotent derivation $\delta=\partial/\partial x$ on a polynomial ring 
$B=k[x]$. Taking $1/x$ as an element $\xi$, it follows that $\sum_{g\in G_a(k)}~k\cdot{}^g\xi=\sum_{n \ge 0}k\cdot x^{-n}$, 
which is not of finite dimension. 

We can extend the notion of algebraicity to regular vector fields 
on an affine algebraic variety. 
Namely, let $D$ be a $k$-derivation on an affine domain $B$. 
We say that $D$ is an {\em algebraic derivation} of $B$ if for 
every element $b \in B$, 
the $k$-vector subspace $\sum_{n\ge 0}k\cdot D^n(b)$ has finite dimension. 
A typical example of such an algebraic derivation 
which is not locally nilpotent is the Euler derivation 
$D=\sum_{i=1}^n x_i\frac{\partial}{\partial x_i}$ 
on a polynomial ring $B=k[x_1,\ldots,x_n]$ and a composite of a derivation 
of Euler type and a locally nilpotent derivation. 

Our objective in this article is to show that an algebraic derivation on an affine $k$-domain is something very close 
to the derivation of the above type. Furthermore, we try to prove some structure theorems on the affine domain 
provided $B$ has an algebraic derivation. We first show that an algebraic derivation $D$ on an affine domain $B$ gives 
a graded ring structure on $B$ (see Theorem \ref{Theorem 1.5}). Then by making use of the graded ring structure and various 
properties of the algebraic derivations, we characterize polynomial rings of dimension two (see Theorems \ref{Theorem 2.4}, 
\ref{Theorem 2.5} and \ref{Theorem 2.6}). In section three, it is shown that the Euler derivation 
$D=x\frac{\partial}{\partial x}+y\frac{\partial}{\partial y}$ on a polynomial ring $k[x,y]$ induces an algebraic 
derivation on $k[x,y]^G$ for every finite subgroup of $\GL(2,k)$. 

For an integral domain $B$, the quotient field of $B$ is denoted by $Q(B)$ and the multiplicative group consisting of 
invertible elements of $B$ by $B^\ast$. For a $k$-derivation $D$ on a $k$-domain $B$, we denote by $\wt{D}$ the extension 
of $D$ to $Q(B)$. We denote by $\Ker D$ and $\Ker \wt{D}$ the kernels of $D$ and $\wt{D}$ respectively, 
which are the subring of $B$ and the subfield of $Q(B)$.

\section{Algebraic derivations}

Let $k$ be an algebraically closed field of characteristic $0$ and 
let $B$ be an affine domain over $k$. 
Examples of algebraic derivations on $B$ are 

\begin{enumerate}
\item[(1)] a locally nilpotent derivation,
\item[(2)] a derivation of Euler type on a polynomial ring 
$B=k[x_1,\ldots,x_n]$.
\[
D=\lambda_1 x_1\frac{\partial}{\partial x_1}+\cdots+\lambda_n x_n\frac{\partial}{\partial x_n}, \quad 
\lambda_1,\ldots,\lambda_n \in k 
\]
\end{enumerate}
Our first purpose is to give a structure theorem of an affine domain $B$ with a non-trivial algebraic derivation $D$.  
Given a derivation $D$, we consider a $k$-algebra homomorphism $\varphi_D : B \to B[[t]]$ defined by 
\[
\varphi_D(b)=\sum_{n \ge 0}\frac 1{n!}D^n(b)t^n \quad \text{for}~ b \in B.
\]
Note that the homomorphism $\varphi_D$ can be defined for any linear map $D : B \to B$. 

For $\lambda \in k$, let $\Delta_\lambda=D-\lambda\cdot 1_B$, where $1_B$ signifies the identity morphism of $B$. 
We often write $\Delta_\lambda=D-\lambda$. 

\begin{lem}\label{Lemma 1.1}
For every $\lambda \in k$, 
we have $\varphi_D =e^{\lambda t}\varphi_{\Delta_\lambda}$. 
\end{lem}

\Proof
Set $\Delta=\Delta_\lambda$ for simplicity. Then $D=\Delta+\lambda\cdot 1_B$. For $b\in B$, 
we compute as follows: 
\begin{eqnarray*}
\lefteqn{\varphi_D(b) = \sum_{n \ge 0}\frac{1}{n!}~(\Delta+\lambda)^n(b)~t^n} \\
&& =b\big(1+\lambda t+\frac{1}{2!}\lambda^2 t^2+\cdots+\frac{1}{n!}\lambda^n t^n+\cdots\big) \\
&& \quad +\Delta(b)t \big(1+\lambda t+\cdots+\frac{1}{(n-1)!}\lambda^{n-1} t^{n-1}+\cdots\big) + \cdots \\
&& \quad +\frac{1}{m!}\Delta^m(b)t^m\big(1+\lambda t+\cdots+\frac{1}{(n-m)!}\lambda^{n-m} t^{n-m}+\cdots\big) + \cdots \\
&&= e^{\lambda t}\big\{b+\Delta(b)t+\frac 1{2!}\Delta^2(b)t^2+\cdots+\frac{1}{m!}\Delta^m(b)t^m+\cdots \big\} \\
&&= e^{\lambda t}\varphi_\Delta(b). 
\end{eqnarray*}
\QED

For $\lambda \in k$, set 
\[
B_\lambda=\{b \in B\mid (D-\lambda)^n(b)=0 ~~\text{for some positive integer}~~ n\}.
\]
Then it is clear that $b \in B_\lambda$ if and only if $\varphi_{\Delta_\lambda}(b) \in B[t]$. In particular,  
\[
B_0=\{b\in B \mid \varphi_D(b)\in B[t]\}.
\]

By making an essential use of Lemma \ref{Lemma 1.1}, we prove the following lemmas.

\begin{lem}\label{Lemma 1.2}
We have the following assertions. 
\begin{enumerate}
\item[(1)] $B_0$ is a $k$-subalgebra of $B$ and 
$D$ is locally nilpotent on $B_0$. 
\item[(2)] For any $\lambda \in k$, $B_\lambda$ is a $B_0$-module. 
\end{enumerate}
\end{lem}
\Proof
(1) The assertion follows from the fact that $\varphi_D$ is a $k$-algebra homomorphism. 

(2) Let $a \in B_0$ and $b \in B_\lambda$. Then 
$$\varphi_D(ab)=\varphi_D(a)\varphi_D(b)
=e^{\lambda t}\varphi_D(a)\varphi_{\Delta_\lambda}(b).$$ 
On the other hand, we have 
$\varphi_D(ab)=e^{\lambda t}\varphi_{\Delta_\lambda}(ab)$. 
Hence it follows that 
$\varphi_{\Delta_\lambda}(ab)
=\varphi_D(a)\varphi_{\Delta_\lambda}(b)\in B[t]$, 
which implies that $\Delta_\lambda^n(ab)=0$ for some $n>0$. 
\QED

\begin{lem}\label{Lemma 1.3}
For $\lambda, \mu \in k$, $B_\lambda B_\mu \subseteq B_{\lambda+\mu}$. 
\end{lem}
\Proof
Let $b \in B_\lambda$ and $c \in B_\mu$. Then 
\[
\varphi_D(bc) = \varphi_D(b)\varphi_D(c) = e^{(\lambda+\mu)t}\varphi_{\Delta_\lambda}(b)\varphi_{\Delta_\mu}(c)
\]
and 
$$\varphi_D(bc)=e^{(\lambda+\mu)t}\varphi_{\Delta_{\lambda+\mu}}(bc).$$
Hence $\varphi_{\Delta_{\lambda+\mu}}(bc)
=\varphi_{\Delta_\lambda}(b)\varphi_{\Delta_\mu}(c)\in B[t]$. 
This implies that 
\[
(D-\lambda-\mu)^n(bc)=0\]
for a sufficiently large $n$, and hence $bc \in B_{\lambda+\mu}$. 
\QED

\begin{lem}\label{Lemma 1.4}
If $\lambda \neq \mu$, then $B_\lambda \cap B_\mu =0$. 
\end{lem}
\Proof
Suppose that $k$ is a subfield of the complex number field  $\C$ . Let $b \in B_\lambda \cap B_\mu$. 
Then $\varphi_D(b)=e^{\lambda t}\varphi_{\Delta_\lambda}(b)=e^{\mu t}\varphi_{\Delta_\mu}(b)$. 
Hence we have 
\[
e^{(\lambda -\mu)t}=\frac{\varphi_{\Delta_\mu}(b)}{\varphi_{\Delta_\lambda}(b)},
\]
where $\varphi_{\Delta_\lambda}(b), \varphi_{\Delta_\mu}(b) \in B[t]$. Hence the function in $t$ on the right 
hand side is a rational function on $\C$. Meanwhile, the function on the left hand side is an entire function 
which does not have zeros and poles. Hence it follows that $\varphi_{\Delta_\mu}(b)$ and $\varphi_{\Delta_\lambda}(b)$ 
are constants, namely, $(D-\lambda)(b)=0$ and $(D-\mu)(b)=0$. So, $D(b)=\lambda b=\mu b$, which implies $b=0$. 

We may assume that $k$ is a subfield of $\C$. Since $B$ is finitely generated over $k$, $B$ is isomorphic to 
the residue ring $k[x_1,\ldots,x_n]/(f_1,\ldots,f_m)$. Let $\alpha_1,\ldots,\alpha_r$ be the coefficients appearing 
in polynomials $f_1,\ldots,f_m \in k[x_1,\ldots,x_n]$, and let $k_0$ be the field $\Q(\alpha_1,\ldots,\alpha_r)$. 
Let 
\[
C=k_0[x_1,\ldots,x_n]/(f_1,\ldots,f_m).
\]
Then $B=C\otimes_{k_0}k$. Adjoining to $k_0$ the elements $\lambda, \mu$ and the coefficients of $D(\ol{x}_i)$ 
for $1 \le i \le n$ when they are expressed in the forms of polynomials in $\ol{x}_1,\ldots,\ol{x}_n$, where 
$\ol{x}_i$ is the residue class of $x_i$ in $B$, we may assume that $\varphi_{\Delta_\lambda}(b)$ and 
$\varphi_{\Delta_\mu}(b)$ are polynomials in $t$ with coefficients in $k_0$. Since $k_0$ is finitely generated over 
$\Q$, we may embed $k_0$ into $\C$ and replace $k$ by the algebraic closure of the embedded $k_0$ in $\C$. 
Thus we may assume that $k \subseteq \C$. 
\QED

For an element $b \in B_\mu$, we define the {\em $\mu$-height} of $b$ as the non-negative integer $r$ 
such that $(D-\mu)^r(b) \ne 0$ and $(D-\mu)^{r+1}(b)=0$. Note that $\mu$-height is defined only for elements 
of $B_\mu$. The first structure theorem on an affine domain with an algebraic derivation is stated as follows. 

\begin{thm}\label{Theorem 1.5}
Let $D$ be an algebraic derivation on an affine domain $B$. Then $B=\bigoplus_{\lambda \in k}B_\lambda$, 
which is a graded ring over $B_0$. Furthermore, $D(B_\lambda) \subseteq B_\lambda$.
\end{thm}

\Proof
Let $b\in B$. Then the vector space $V=\sum_{n\ge 0}k D^n(b)$ has finite dimension. 
We may choose an integer $n > 0$ so that 
\[
\{b,D(b), D^2(b),\ldots, D^{n-1}(b)\}
\]
is a $k$-basis of $V$ and $D^n(b)$ is expressed by a linear combination of $b, D(b), \ldots, D^{n-1}(b)$.
Then $D$ is a $k$-linear endomorphism on $V$ and $V$ decomposes into a direct sum $V=\oplus_\lambda V_\lambda$ 
where $V_\lambda=\{v \in V \mid (D-\lambda)^n(v)=0, n \gg 0\}$. Since $V_\lambda \subseteq B_\lambda$, 
it follows that $b\in \bigoplus_\lambda B_\lambda$. The previous lemmas show that $B=\bigoplus_\lambda B_\lambda$ 
is a graded ring over $B_0$. As for the second assertion, let $b \in B_\lambda$. Let $r$ be the $\lambda$-height of $b$. 
If $r=0$, then $D(b)=\lambda b \in B_\lambda$. Suppose that $r > 0$. Since $(D-\lambda)^{r+1}(b)=0$, the element 
$(D-\lambda)(b) \in B_\lambda$ and the $\lambda$-height of $(D-\lambda)(b)$ is $r-1$. By induction, we may assume 
$(D-\lambda)(b) \in B_\lambda$. Then $D(b)=\lambda b+(D-\lambda)(b) \in B_\lambda$. So, $D(B_\lambda) \subseteq 
B_\lambda$. 
\QED

We look into the properties of the subring $B_0$.

\begin{lem}\label{Lemma 1.6}
Let $D$ be an algebraic derivation on a normal affine domain $B$. Then $B_0$ is integrally closed. 
\end{lem}
\Proof
Suppose that $\xi \in Q(B_0)$ is integral over $B_0$. Then it follows that $\xi \in B$ since $B$ is normal. 
The element $\xi \in B$ satisfies 
\[
\xi^n+\alpha_1\xi^{n-1}+\cdots +\alpha_n=0,
\]
where each $\alpha_i$ is an element of $B_0$. Hence it follows that 
\[
\varphi_D(\xi)^n+\varphi_D(\alpha_1)\varphi_D(\xi)^{n-1}+\cdots +\varphi_D(\alpha_n)=0.
\]
So, $\varphi_D(\xi)$ is integral over $k[\varphi_D(\alpha_1),\ldots, \varphi_D(\alpha_n)]\subset B[t]$. 
Since $B[t]$ is integrally closed and $\varphi_D(\xi) \in Q(B[t])$, it follows that $\varphi_D(\xi)\in B[t]$. 
Hence $\xi \in B_0$. 
\QED

We observe the following two examples. One deals with the Euler derivation and the other does a composite of 
the Euler derivation and a locally nilpotent derivation. 

\begin{example}\label{Example 1.7}{\em 
Let $B=k[x_1,\ldots,x_n]$ be a polynomial ring in $n$ variables and let $D=\sum_{i=1}^n x_i\frac{\partial}{\partial x_i}$ 
be the Euler derivation on $B$. Then $D$ is algebraic and $B=\bigoplus_{d\ge 0}B_d$, where $B_d$ is the $k$-vector space 
spanned by monomials $x_1^{m_1}\cdots x_n^{m_n}$ of total degree $d$. In particular, $B_0=k$. Meanwhile, $\Ker \wt{D}$ 
is $k\big(\frac{x_2}{x_1},\ldots,\frac{x_n}{x_1}\big)$ which has transcendence degree $n-1$ over $k$.
\QED}
\end{example}

\begin{example}\label{Example 1.8}{\em 
Let $B=k[x,y]$ and let $D=x\frac{\partial}{\partial x}+\frac{\partial}{\partial y}$. For non-negative integers 
$m$ and $n$, we have 
\begin{equation}
D(x^ny^m)=nx^ny^m+mx^ny^{m-1}. \label{1}
\end{equation}
Hence for any non-negative integer $r$, we have 
$D^r(x^ny^m)\in k\cdot x^ny^m+k \cdot x^ny^{m-1}+\cdots +kx^n$. 
Thus it follows that the derivation $D$ on $B$ is algebraic.   
By (\ref{1}), we have $(D-n)(x^ny^m)=mx^ny^{m-1}$ and 
\[
(D-n)^{m+1}(x^ny^m)=0.
\]
Hence if $f(y)\in k[y]$ with $\deg f(y)=m$, then 
\[
(D-n)^{m+1}(x^nf(y))=0
\]
and $x^nf(y) \in B_n$. Since any element $b \in B$ is written as $b=f_0(y)+xf_1(y)+\cdots +x^nf_n(y)$ where $f_i(y)\in k[y]$ 
for $1 \le i \le n$, it follows that $B=\bigoplus_{n \in \Z_{\ge 0}}B_n$ and $B_n=x^nk[y]$. 
The ring $B$ is graded by the monoid $\Z_{\ge 0}$ of non-negative integers. \QED}
\end{example}

Let $D$ be an algebraic derivation on a normal affine domain $B$ over $k$. 
Let $\Lambda=\{\lambda \in k\mid B_\lambda \neq 0\}$, which is a monoid under the addition of $(k,+)$, 
and let $M$ be the abelian subgroup of $(k,+)$ generated by $\Lambda$. We call $\Lambda$ (resp. $M$) the 
{\em monoid} (resp. {\em abelian group}) associated to $D$.

\begin{thm}\label{Theorem 1.9}
Under the notation and assumption as above, 
suppose that $M$ is a totally ordered abelian group with ordering $<$ 
in the sense that $\lambda \ge 0$ for every $\lambda \in \Lambda$ and 
if $\lambda < \mu$ then $\lambda +\nu <\mu +\nu$ for any $\nu \in M$. 
\footnote{The condition that if $\lambda < \mu$ then 
$\lambda +\nu <\mu +\nu$ for any $\nu \in M$ 
follows from that $M$ is a finitely generated subgroup of the additive group $k$. }
Then the following assertions hold. 
\begin{enumerate}
\item[(1)] $Q(B_0)\cap B=B_0$. 
\item[(2)] $Q(B_0)$ is algebraically closed in $Q(B)$. 
If $B_0 \subsetneqq B$, then we have $\transdeg_kQ(B_0) < \transdeg_kQ(B)$. 
\item[(3)] $B_0$ is factorially closed in $B$. 
Hence if $B$ is factorial, then so is $B_0$. 
\item[(4)] Suppose that $B$ is a polynomial ring of $\dim \le 3$. Then $B_0$ is either a polynomial ring 
or $k$. 
\end{enumerate}
\end{thm}
\Proof
(1)\ \ Suppose that $b \in Q(B_0)\cap B$. Then $b$ is written as $a_0b=a_1$ with $a_0, a_1 \in B_0$. 
By Theorem \ref{Theorem 1.5}, $b=b_{\lambda_1}+\cdots +b_{\lambda_n}$ 
with $b_{\lambda_i}\in B_{\lambda_i}$ for $1 \le i \le n$. 
Since $a_0b=a_0b_{\lambda_1}+\cdots +a_0b_{\lambda_n}=a_1$, it follows that 
$\lambda_i=0$ for some $i$ and $b=b_{\lambda_i}\in B_0$. 

(2)\ \ Suppose that $\xi \in Q(B)$ is algebraic over $Q(B_0)$. 
Then there exists $a_0 \in B_0$ such that $a_0\xi$ is integral over $B_0$. 
Since $B_0 \subset B$ and $B$ is integrally closed by assumption, 
$a_0\xi \in B$. The element $a_0\xi$ satisfies 
\begin{equation}
(a_0\xi)^n+c_1(a_0\xi)^{n-1}+\cdots +c_n=0 \label{2} 
\end{equation}
where $c_i\in B_0$ for $1 \le i \le n$. Write 
\[
a_0\xi=b_{\lambda_1}+\cdots + b_{\lambda_m}, \quad b_{\lambda_1}\neq 0,~
b_{\lambda_m}\neq 0, \quad \lambda_1 < \cdots <\lambda_m.
\]
If $\lambda_m >0$, then the term $b_{\lambda_m}^n$, 
which has the highest order, cannot be cancelled by the other terms
in the equation (\ref{2}). Hence $\lambda_m \le 0$. 
Similarly, $\lambda_1\ge 0$. So, $a_0\xi \in B_0$ and $\xi \in Q(B_0)$.  
Suppose that $B_0 \subsetneqq B$. 
Then there exists an element $b \in B_\lambda$ with $\lambda \ne 0$. 
Then $b$ is algebraically independent over $Q(B_0)$. 
For otherwise, $b$ is algebraic over $Q(B_0)$ and hence $b \in Q(B_0)$. 
Then $b \in Q(B_0)\cap B=B_0$, which contradicts the choice of $b$.   

(3)\ \ Suppose that $a=b_1b_2\in B_0$ with non-zero $b_1, b_2 \in B$. 
Write $b_1=b_{\lambda_1}+\cdots +b_{\lambda_r}$ and 
$b_2=b_{\mu_1}+\cdots +b_{\mu_s}$ where 
$b_{\lambda_i}\in B_{\lambda_i}$ and $b_{\mu_j}\in B_{\mu_j}$ 
with $\lambda_1 < \cdots < \lambda_r$ and $\mu_1<\cdots <\mu_s$. 
We may assume that $b_{\lambda_1}b_{\lambda_r}b_{\mu_1}b_{\mu_s}\neq 0$. 
Then the highest term of $b_1b_2$ is $b_{\lambda_r}b_{\mu_s}$ and 
the lowest is $b_{\lambda_1}b_{\mu_1}$. 
Hence $b_{\lambda_r}b_{\mu_s}=b_{\lambda_1}b_{\mu_1}=a$ and so, $b_1=b_{\lambda_1}$ and $b_2=b_{\mu_1}$. 
Since $\lambda_1\ge 0$, $\mu_1\ge 0$ and $\lambda_1+\mu_1=0$, 
it follows that $b_1, b_2 \in B_0$ 
and $B_0$ is factorially closed in $B$. 

(4)\ \ We may assume that $B_0 \subsetneqq B$ and $B_0 \ne k$. Then $\transdeg_kQ(B_0)$ $< \transdeg_kQ(B) \le 3$. 
Since $\dim B \le 3$, it follows by the assertion (1) and a result of Zariski \cite{Z} that $B_0$ is an affine domain. 
Let $X=\Spec B$ and $Y=\Spec B_0$. Then the inclusion $B_0\hookrightarrow B$ induces a dominant morphism $p : X \to Y$. 
Since $B$ is factorial and $B^*=k^*$, $B_0$ is factorial by the assertion (3) and $B_0^*=k^*$. 
If $\dim B_0=1$, then $B_0$ is a polynomial ring. Suppose that $\dim B_0=2$. Then $p$ is equi-dimensional over $p(X)$. 
In fact, suppose that there is an irreducible fiber component of dimension $2$ in $p$. Then there exists a prime element 
$b \in B$ such that $bB\cap B_0= \fm$ is a maximal ideal. Then $bb_1=a \in B_0$. Since $B_0$ is factorially closed in $B$, 
it follows that $b \in B_0$. Hence $bB\cap B_0=bB_0$, which is a contradiction. Then it follows that $Y$ is isomorphic 
to $\A^2$ or a Platonic fiber space $\A^2/G$ (see \cite[Theorem 3]{MM}). If $G$ is nontrivial, $p$ splits to 
$X \to \A^2 \to \A^2/G=Y$. In fact, $Y$ has a unique singular point, say $P$ and the smooth locus $Y^\circ$ has 
the universal covering $\A^2\setminus\{O\}$, where $O$ is the point of origin. Since $p^{-1}(P)$ is either empty 
or of dimension one, $\A^3\setminus p^{-1}(O)$ is simply connected. Hence the fiber product 
\[
Z=(\A^3\setminus p^{-1}(P))\times_{Y^\circ}(\A^2\setminus\{O\})
\]
splits into a disjoint union of copies of $\A^3\setminus p^{-1}(P)$. Then the restriction of the second projection 
$p_2 : Z \to \A^2\setminus\{O\}$ to a connected component of $Z$ provides the above splitting of $p$. 
This contradicts the fact that $Q(B_0)$ is algebraically closed in $Q(B)$. 
\QED

\begin{example}\label{Example 1.10}{\em
Let $B=k[x,y]$ and 
\[
D=x^{n-1}y^n(x\frac{\partial}{\partial x}-y\frac{\partial}{\partial y}) \quad (n \ge 1).
\]
We have by computation 
\begin{eqnarray*}
D(x) &=& x^ny^n, ~~ D^2(x)=0 \quad \text{and} \\
D^k(y)&=&(-1)^kk!x^{k(n-1)}y^{kn+1} \quad (k \ge 1)
\end{eqnarray*}
Thus $\sum_{j=0}^\infty k\cdot D^j(y)$ is infinite-dimensional, and $D$ is not algebraic.

Next, compute $D^\ell(x^iy^j)$. For $0 \le i < j$, we have 
\[
D^\ell(x^iy^j)=(i-j)(i-j-1)\cdots(i-j-\ell+1)x^{i+\ell(n-1)}y^{j+\ell n}.
\]
For $i=j$, we have $D(x^iy^i)=0$. 
For $i > j$, it follows that $x^iy^j=x^{i-j}(xy)^j \in B_0$. 
Hence $B_0=k[x,xy]$. Therefore, we have $Q(B_0)=Q(B)$ and 
$$Q(B_0)\cap B=Q(B)\cap B=B \supsetneqq B_0.$$
If we write $z=xy$, then $B_0=k[x,z]$ and
$D|_{B_0}=z^n\frac{\partial}{\partial x}$. 
\QED}
\end{example}

We can raise the following natural question. 
\svskip

\noindent
{\bf Question.}\ \ Let $B$ be an affine domain over $k$ 
with an algebraic derivation $D$. Is $B_0$ finitely generated over $k$?
\svskip

When $D$ is locally nilpotent, it is clear that $B_0=B$ is finitely generated. 
So, we are interested in the case that 
$D$ is not a locally nilpotent derivation. 
If the monoid $\Lambda$ is $\Z_{\ge 0}$, 
then $B_0$ is an affine domain and every $B_\mu$ is a finitely 
generated $B_0$-module. 
\svskip

An element $\lambda$ of the monoid $\Lambda$ is said to be {\em primitive} if $\lambda=\mu+\nu$ with 
$\mu, \nu \in \Lambda$ then either $\mu=0$ or $\nu=0$. 

\begin{lem}\label{Lemma 1.11}
Let $\lambda$ be a primitive element of the associated monoid $\Lambda$ and let $\xi$ be an element of 
$B_\lambda$ such that $(D-\lambda)(\xi)=0$. Suppose that the associated abelian group $M$ is totally 
ordered, $B$ is a factorial domain with $B^\ast=k^\ast$ and $\Ker D=k$. Then $\xi$ is irreducible in $B$.
\end{lem}
\Proof
Suppose that $\xi$ is reducible. Write $\xi=b_1b_2$ with $b_1, b_2 \in B$ and $\gcd(b_1,b_2)=1$. 
Since $\xi$ is a homogeneous element, both $b_1$ and $b_2$ are homogeneous. Since $\lambda$ is primitive 
in $\Lambda$, we may assume that $b_1 \in B_0$ and $b_2 \in B_\lambda$. Since $(D-\lambda)(b_1b_2)=0$, 
we have 
\[
b_1(D(b_2)-\lambda b_2)=-b_2D(b_1).
\]
Since $\gcd(b_1,b_2)=1$, $D(b_1)$ is divisible by $b_1$. By Theorem \ref{Theorem 1.9}, $B_0$ is 
a factorial affine domain over $k$. Since $D(B_0) \subseteq B_0$ and $D$ is locally nilpotent on $B_0$, 
it follows that $b_1 \in \Ker D=k$ (see \cite{F}). This is a contradiction.
\QED

\section{Algebraic derivations on $k[x,y]$}

Let $B$ be an affine domain with an algebraic derivation $D$. Then $D$ extends uniquely to a derivation 
$\tilde D$ on $K=Q(B)$. In the present section, we intend to characterize a two-dimensional polynomial 
ring $k[x,y]$ and algebraic derivations on $k[x,y]$ in terms of the subring $B_0$ and 
the associated monoid $\Lambda$, though our results are partial. We begin with the following result. 

\begin{lem}\label{Lemma 2.1}
Let $R$ be an affine domain of dimension one defined over a field $K_0$ of characteristic zero 
which is not necessarily algebraically closed. Assume that $K_0$ is algebraically closed in $Q(R)$ 
and that $R$ is a graded ring $R=\bigoplus_{n \ge 0}R_n$ with $R_0=K_0$. Then the following assertions 
hold.
\begin{enumerate}
\item[(1)] 
$R$ is a polynomial ring $K_0[\xi]$ such that $R_n=K_0\xi^n$ 
for every $n \ge 0$. 
\item[(2)]
Suppose that the graded ring structure on $R$ is given by 
an algebraic $K_0$-derivation $D$ and 
the associated monoid $\Lambda$ is $\Z_{\ge 0}\lambda$ 
for some $\lambda \in \Lambda$. 
Then for every $\mu \in \Lambda$, we have 
\[
R_\mu=\{\eta \in R \mid (D-\mu)(\eta)=0\}.
\]
\end{enumerate}
\end{lem}

\Proof
(1)\ \ Let $\xi$ be a nonzero element of $R_1$ and let $A=K_0[\xi]$. Then $\xi$ is algebraically 
independent over $K_0$ and hence $A$ is a graded subalgebra of dimension one in $R$. Since $\dim R=1$, 
the field extension $Q(R) \supset Q(A)$ is algebraic. Let $\eta \in R_n$. Then there exist elements 
$f_0(\xi),f_1(\xi),\ldots,f_r(\xi) \in A$ with $\gcd(f_0(\xi),\ldots,f_r(\xi))=1$ such that
\begin{equation}
f_0(\xi)\eta^r+f_1(\xi)\eta^{r-1}+ \cdots+f_r(\xi)=0   \label{3}
\end{equation} 
yields a minimal equation of $\eta$ over $Q(A)$. Write 
\[
f_0(\xi)=\xi^m+\text{(terms of lower degree in $\xi$)}.
\]
Consider the homogeneous part of degree $(m+rn)$ in the equation (\ref{3}) 
and obtain
\begin{equation*}
\xi^m\eta^r+c_1\xi^{m+n}\eta^{r-1}+\cdots+c_r\xi^{m+nr}=0, 
\quad c_1, \ldots,c_r \in K_0. 
\end{equation*}
Hence we have 
\begin{equation}
\eta^r+c_1\xi^{n}\eta^{r-1}+\cdots+c_r\xi^{nr}=0.  \label{4}
\end{equation}
By the minimality of the equation (\ref{3}) over $Q(A)$, 
the equation (\ref{3}) coincides with the equation (\ref{4}) up to $K_0^*$. 
Hence $m=0$ and the equation (\ref{3}) is written as 
\begin{equation}
\eta^r+c_1\xi^n\eta^{r-1}+\cdots+c_r\xi^{nr}=0. \label{5}
\end{equation}
The equation (\ref{5}) yields an algebraic equation
\[
\left(\frac{\eta}{\xi^n}\right)^r+c_1\left(\frac{\eta}{\xi^n}\right)^{r-1}+\cdots +c_r=0.
\]
Since $\frac{\eta}{\xi^n}$ is an element of $Q(R)$ and $K_0$ is algebraically closed in $Q(R)$, 
it follows that $\frac{\eta}{\xi^n}$ is an element of $K_0$. This implies that $\eta \in K_0\xi^n$, 
and hence $R=A$. 

(2)\ \ Note that $\lambda$ is a minimal element of $\Lambda$ other than $0$. Let $\xi$ be a nonzero 
element of $R_\lambda$ such that $D(\xi)=\lambda\xi$. Then $R$ is a polynomial ring $K_0[\xi]$ by 
the assertion (1). Let $\mu \in \Lambda$ and write $\mu=n\lambda$. Then every element $\eta$ of $R_\mu$ 
is of the form $c\xi^n$ with $c \in K_0$.  Hence we compute
\[
D(\eta)=D(c\xi^n)=cn\xi^{n-1}D(\xi)=cn\lambda \xi^n=(n\lambda)\eta=\mu\eta.
\]
Thence follows the assertion (2).
\QED

The following result determines an algebraic derivation on $B=k[x]$.

\begin{prop}\label{Proposition 2.2}
Let $B=k[x]$ be a polynomial ring of dimension one and 
let $D$ be a nontrivial algebraic derivation on $B$. 
Then, after a suitable change of variable, 
either $D=c\frac{d}{dx}$ or $D=cx\frac{d}{dx}$, where $c \in k^\ast$.
\end{prop}
\Proof
If $B=B_0$ then $D$ is locally nilpotent, 
and $D=c\frac{d}{dx}$ with $c \in k^\ast$. 
Suppose that $B \ne B_0$. Let $\xi \in B\backslash B_0$ 
such that $D(\xi)=\lambda\xi$ with $\lambda \in k^\ast$, 
and write $\xi=f(x) \in k[x]$. 
Then we have $D(\xi)=f'(x)D(x)=\lambda f(x)$. 
Since $\deg f(x) > 0$, it follows that $D(x)$ is a linear polynomial in $x$. Write $D(x)=cx+d$ with $c \ne 0$. 
Then $D(cx+d)=c(cx+d)$. By replacing $x$ with $cx+d$, we may assume that $D(x)=cx$ with $c \in k^\ast$. 
Then $B=\bigoplus_{n \ge 0}B_n$ with $B_n=kx^n$ is the homogenoeous decomposition of $B$ with respect to $D$. 
Hence $D=cx\frac{d}{dx}$.
\QED

In the case $\dim B=2$, if we assume that the monoid $\Lambda$ 
associated to an algebraic derivation $D$ is isomorphic to $\Z_{\ge 0}$, 
we have the following result. 

\begin{thm}\label{Theorem 2.3}
Let $B$ be an affine domain of dimension two. 
Then $B$ is isomorphic to $k[x,y]$ if and only if the following 
conditions are satisfied.
\begin{enumerate}
\item[(1)]
$B$ is a factorial domain with $B^\ast=k^\ast$.
\item[(2)]
$B$ has a nontrivial algebraic $k$-derivation $D$ such that the abelian group $M$ generated by the associated 
monoid $\Lambda$ is a totally ordered abelian group, $\dim B_0 \ge 1$ and $\Lambda$ is isomorphic to $\Z_{\ge 0}$ 
provided $\dim B_0=1$.
\end{enumerate}
\end{thm}
\Proof
Since the \lq\lq only if'' part is easy, we prove only the \lq\lq if'' part. 
If $\dim B_0=2$ then $B=B_0$. 
In fact, if $B \supsetneqq B_0$, take any element $\xi \in B\backslash B_0$. 
Then $\xi$ is algebraic over $Q(B_0)$. 
Since $Q(B_0)$ is algebraically closed in $Q(B)$ and $Q(B_0)\cap B=B_0$ by 
Theorem \ref{Theorem 1.9}, it follows that $\xi \in B_0$, a contradiction. Hence $B=B_0$. Then $D$ 
is locally nilpotent on $B$. By an algebraic characterization of $k[x,y]$, $B$ is a polynomial ring $k[x,y]$
\cite[Theorem 2.2.1]{M1}. Suppose that $\dim B_0=1$. Then $B_0$ is an affine domain of dimension one such that 
$B_0$ is factorial and $B_0^\ast=k^\ast$. Hence $B_0=k[x]$. By the assumption, $B=\bigoplus_{n \ge 0}B_{n\lambda}$. 
Let $K_0=Q(B_0)$ and $R=B\otimes_{B_0}K_0$. Then $R$ is an affine domain of dimension one over $K_0$ which has 
the graded ring structure $R=\bigoplus_{n\ge 0}R_n$ with $R_0=K_0$, where $R_n=B_{n\lambda}\otimes_{B_0}K_0$. 
By Lemma \ref{Lemma 2.1}, $R$ is a polynomial ring $K_0[\xi]$. Hence the affine surface $\Spec B$ has an 
$\A^1$-fibration \cite{GKM}, and $B$ is a polynomial ring in two variables by the algebraic characterization 
of $k[x,y]$ \cite{M1}. 
\QED  

Under some additional conditions, we can further determine 
an algebraic derivation. 

\begin{thm}\label{Theorem 2.4}
Let $B$ be an affine domain of dimension two with an algebraic derivation $D$. 
Then $B=k[x,y]$ with 
\[
D=\frac{\partial}{\partial x}+\lambda y\frac{\partial}{\partial y}
\] 
for $\lambda \in k^*$ if and only if the following conditions are satisfied. 
\begin{enumerate}
\item[(1)] $B$ is factorial and $B^*=k^*$. 
\item[(2)] $\Ker \wt{D}=\Ker D=k$ and $B_0 \neq k$. 
\item[(3)] The monoid $\Lambda$ associated to $D$ is the set $\Z_{\ge 0}\lambda$ of non-negative multiples of $\lambda \in k^\ast$. 
\end{enumerate}
\end{thm}
\Proof ``Only if'' part. 
Suppose that $B=k[x,y]$ and 
$D=\frac{\partial}{\partial x}+\lambda y\frac{\partial}{\partial y}$. 
Then we can compute for every $\ell \ge 1$
\[
D^\ell(x^my^n)= n^\ell\lambda^\ell x^my^n+\sum_{0\le r < m}c_rx^ry^n, \quad c_r \in k.
\]
Hence $x^my^n \in B_0$ if and only if $n=0$. This implies that $B_0=k[x]$. 
Further, $B_\mu\neq 0$ if and only if $\mu \in \Z_{\ge 0}\lambda$ and $B_{n\lambda}=B_0y^n$ for a positive integer $n$. 
In order to find elements in $\Ker \wt{D}$, write an element $\xi \in Q(B)$ as $\xi=\frac{f}{g}$ with 
$f, g \in k[x,y]$ and $\gcd(f,g)=1$. Then $\wt{D}(\xi)=0$ implies $f \mid D(f)$ and $g \mid D(g)$. 
Write $f=a_0+a_1y+\cdots+a_my^m$ with the $a_i \in k[x]$. Then it is straightforward to show that $f \mid D(f)$ only if 
$f$ is a monomial in $y$ with coefficients in $k$. Similarly, $g$ is a monomial in $y$. Then $\wt{D}(\xi)=0$ 
if and only if $\xi \in k$. The rest of the assertion is easy to show. 

``If'' part. By Theorem \ref{Theorem 2.3}, $B$ is a polynomial ring $k[x,y]$. 
If $B=B_0$, then $D$ is locally nilpotent and 
$D=f(y)\frac{\partial}{\partial x}$ for a suitable choice of coordinates. 
It then follows that $\Ker D\neq k$, which is a contradiction. 
Hence $B\neq B_0$ and $\dim B_0=1$ by Theorem \ref{Theorem 1.9}. 
Since $B_0$ is a factorial affine domain of dimension one with $B_0^*=k^*$, 
$B_0$ is a polynomial ring and we may put $B_0=k[x]$. 
We may also assume that $D(x)=1$. For every $\mu \in \Lambda$, 
there exists $\xi_\mu \in B_\mu$ such that $D(\xi_\mu)=\mu \xi_\mu$. 
In fact, take a nonzero element $z \in B_\mu$ and 
let $r$ be the $\mu$-height of $z$. 
Let $\xi_\mu=(D-\mu)^rz$. Then $\xi_\mu\in B_\mu$ satisfies $D(\xi_\mu)=\mu\xi_\mu$. We claim that for every 
$\mu \in \Lambda$
\[
B_\mu=B_0\xi_\mu.
\]
Take any $\eta \in B_\mu$. If $(D-\mu)\eta=0$, then $\wt{D}\big(\frac{\eta}{\xi_\mu}\big)=0$. 
Hence $\frac\eta{\xi_\mu}\in \Ker\wt{D}=k$ and $\eta=c\xi_\mu$ with $c \in k$. 
We can show the claim by induction on the $\mu$-height $r$ of $\eta$. Then $(D-\mu)^r\eta\neq 0$ and $(D-\mu)^{r+1}\eta=0$. 
By the induction hypothesis, there exists a polynomial $f(x)\in k[x]$ such that $(D-\mu)\eta=f(x)\xi_\mu$. 
Let $F(x)\in k[x]$ be a polynomial such that $F'(x)=f(x)$. 
Then 
\[
(D-\mu)(F(x)\xi_\mu)=F'(x)\xi_\mu +F(x)D(\xi_\mu)-\mu F(x)\xi_\mu=f(x)\xi_\mu.
\]
Hence it follows that $(D-\mu)(\eta-F(x)\xi_\mu)=0$, and 
$\eta -F(x)\xi_\mu=c\xi_\mu$ for $c \in k$. 
Thus the claim is verified. We can also verify by induction that $\deg f(x)$ is the $\mu$-height of $\eta$ 
when we write $\eta=f(x)\xi_\mu$.

Since $\Lambda=\Z_{\ge 0}\lambda$ by the hypothesis, it follows that $B_\mu=B_0\xi_\lambda^n$ if $\mu=\lambda n$. 
Hence $B=k[x,\xi_\lambda]$ and $D$ is written as 
$D=\frac{\partial}{\partial x}+\lambda y\frac{\partial}{\partial y}$ with $y=\xi_\lambda$. 
\QED

\smallskip

The following result differs from Theorem \ref{Theorem 2.4} only in the condition (2).
 
\begin{thm}\label{Theorem 2.5}
Let $B$ be an affine domain of dimension two with a nontrivial algebraic derivation $D$. 
Then $B=k[x,y]$ with $D=\lambda y\frac{\partial}{\partial y}$ and $\lambda \in k^*$ 
if and only if the following conditions are satisfied. 
\begin{enumerate}
\item[(1)] $B$ is factorial and $B^*=k^*$. 
\item[(2)] $\Ker \wt{D}=Q(\Ker D) \supsetneqq k$.  
\item[(3)] The monoid $\Lambda$ associated with $D$ is $\Z_{\ge 0}\lambda$ with $\lambda \in k^\ast$. 
\end{enumerate}
\end{thm}
\Proof
``Only if'' part. Clearly, $\Ker D=k[x]=B_0$ and $B=\bigoplus_{n\ge 0}B_0 y^n$ with 
$\Ker \wt{D}=k(x)$ and $\Lambda=\Z_{\ge 0}\lambda$. 

``If'' part. By the conditions (2) and (3), we have the inclusions
\[
k \subsetneqq \Ker D \subseteq B_0 \subsetneqq B. 
\]
Hence $\dim \Ker D=\dim B_0=1$. 
Then every element of $B_0$ is algebraic over $\Ker D$. 
Since $D$ is locally nilpotent over $B_0$, 
it follows that $\Ker D=B_0$ (\cite{F}). 
By Theorem \ref{Theorem 1.9}, $B_0=k[x]$. 
Note that $\lambda$ is a minimal element of $\Lambda$ which is nonzero and $\Lambda=\Z_{\ge 0}\lambda$. 
Let $\xi$ be an element of $B_\lambda$ such that $(D-\lambda)\xi=0$. By dividing $\xi$ by elements of $k[x]$, 
we may assume that $\xi$ is not divisible by any element of $k[x]\setminus k$. 
Let $\eta \in B_\mu$ where $\mu=n\lambda$. Suppose that $(D-\lambda)\eta=0$. 
Then $\wt{D}\big(\frac{\eta}{\xi^n}\big)=0$. Hence $\frac{\eta}{\xi^n}\in Q(B_0)$ and thus we find  
$f(x),g(x)\in B_0$ such that $f(x)\eta=g(x)\xi^n$ and $\gcd(f(x),g(x))=1$. Since $B$ is factorial 
and $\xi$ is not divisible by any element of $k[x]\setminus k$, it follows that $f(x)\in k$ and 
$\eta \in B_0\xi^n$. 

We shall show that $B_\mu=\{\eta \in B \mid (D-\mu)(\eta)=0\}$. Let $K_0=Q(B_0)$ and $R=B\otimes_{B_0}K_0$. 
By Theorem \ref{Theorem 1.9}, $Q(B_0)$ is algebraically closed in $Q(R)$. Then $R$ is an affine 
domain of dimension one over $K_0$ and $D$ extends to a $K_0$-trivial algebraic derivation on $R$ 
because $D(x)=0$. We may assume that $B \subset R$. Let $\eta \in B_\mu$. Then $\eta \in R_\mu$ and 
$(D-\mu)(\eta)=0$ by Lemma \ref{Lemma 2.1}. 

Now we have proved that $R_\mu=B_0\xi^n$. Hence $B=k[x,\xi]$ and $D(\xi)=\lambda\xi$. So, it suffices to take 
$y=\xi$.
\QED

We consider a characterization of an algebraic derivation of Euler type on a polynomial ring 
$k[x_1,\ldots,x_r]$. Let $D$ be an algebraic derivation on an affine domain $B$. 
Let $\lambda_1, \ldots, \lambda_r$ be elements of $k$. We say that $\lambda_1,\ldots,\lambda_r$ are 
{\em numerically independent} if $n_1\lambda_1+\cdots+n_r\lambda_r=0$ with the $n_i \in \Z$ implies 
$n_1=\cdots=n_r=0$. The monoid $\Lambda$ of $D$ is said to be a {\em free} monoid of 
{\em rank} $r$ if there exists numerically independent elements $\lambda_1,\ldots,\lambda_r$ of $k$ such that 
$\Lambda=\sum_{i=1}^r\Z_{\ge 0}\lambda_i$. Then the abelian group $M$ generated by $\Lambda$ is a free abelian 
group $\bigoplus_{i=1}^r\Z\lambda_i$ of rank $r$. Define a lexicographic order in $M$ by setting
\[
m_1\lambda_1+\cdots+m_r\lambda_r < n_1\lambda_1+\cdots+n_r\lambda_r
\]
if there exists $1 \le \ell \le r$ such that $m_i=n_i$ for $1 \le i < \ell$ and $m_\ell <n_\ell$. Then $M$ is 
a totally ordered abelian group with the lexicographic order. 

\begin{thm}\label{Theorem 2.6}
Let $B$ be an affine domain of dimension $r$ with a nontrivial algebraic derivation $D$. 
Then $B=k[x_1,\ldots,x_r]$ and $D=\sum_{i=1}^r\lambda_i x_i\frac{\partial}{\partial x_i}$ with numerically 
independent elements $\lambda_1,\ldots,\lambda_r \in k$ if and only if the following conditions are satisfied.
\begin{enumerate}
\item[(1)]
$B$ is factorial and $B^\ast=k^\ast$.
\item[(2)]
$\Ker D=k$.
\item[(3)]
The monoid $\Lambda$ is a free monoid of rank $r$.
\end{enumerate}
\end{thm}
\Proof
The \lq\lq only if '' part is clear with an algebraic derivation
$D=\sum_{i=1}^r\lambda_i x_i\frac{\partial}{\partial x_i}$, where $\lambda_1, \ldots,\lambda_r$ are 
numerically independent. The monoid $\Lambda$ is then generated by $\lambda_1,\ldots, \lambda_r$ over $\Z_{\ge 0}$ 
and $B_\lambda=kx_1^{m_1}\cdots x_r^{m_r}$ if $\lambda=m_1\lambda_1+\cdots+m_r\lambda_r$. 

We shall prove the \lq\lq if '' part. By the condition (3), the monoid is written as 
$\Lambda=\sum_{i=1}^r \lambda_i\Z_{\ge 0}$ with numerically independent elements $\lambda_1,\ldots, \lambda_r$ 
of $k$. For $1 \le i \le r$, choose elements $\xi_i \in B_{\lambda_i}$ such that $(D-\lambda_i)(\xi_i)=0$. 
Since $\lambda_i$ is primitive in $\Lambda$ and the condition (2) is assumed, $\xi_1,\ldots,\xi_r$ are 
irreducible elements in $B$ by Lemma \ref{Lemma 1.11}. We claim that $\xi_1,\ldots,\xi_r$ are algebraically 
independent over $k$. In fact, let 
\[
\sum_{i_1,\ldots,i_r}c_{i_1\cdots i_r}\xi_1^{i_1}\cdots\xi_r^{i_r}=0, \quad c_{i_1\cdots i_r} \in k
\]
be an algebraic relation of $\xi_1,\ldots,\xi_r$ over $k$. Note that 
\[
\xi_1^{i_1}\cdots\xi_r^{i_r} \in B_{i_1\lambda_1+\cdots+i_r\lambda_r}
\]
and $i_1\lambda_1+\cdots+i_r\lambda_r \ne j_1\lambda_1+\cdots+j_r\lambda_r$ in $\Lambda$ whenever 
$(i_1,\ldots,i_r)\ne (j_1,\ldots,j_r)$. 
Since $B=\bigoplus_{\lambda \in \Lambda}B_\lambda$, it follows that 
\[
c_{i_1\cdots i_r}=0 \quad \text{for}\quad \forall~(i_1,\ldots,i_r).
\]
Let $A=k[\xi_1,\ldots,\xi_r]$. Then $A$ is a graded subalgebra of $B$ of dimension $r$. Hence $Q(B)$ is algebraic 
over $Q(A)$. Let $\eta$ be a homogeneous element of $B_\mu$. Then there exist elements $a_0(\xi), \ldots, a_s(\xi)$ 
of $A$ such that 
\begin{equation}
a_0(\xi)\eta^s+a_1(\xi)\eta^{s-1}+\cdots+a_s(\xi)=0. \label{6}
\end{equation}
We may assume that this gives a minimal equation of $\eta$ over $Q(A)$ 
and $\gcd(a_0, \ldots,a_s)=1$. We can write 
\[
a_0(\xi)=\xi_1^{m_1}\cdots\xi_r^{m_r}+
\text{(lower terms in $\xi$ with the lexicographic order)}.
\]
Write $\mu=\alpha_1\lambda_1+\cdots+\alpha_r\lambda_r$ 
with $\alpha_1, \ldots, \alpha_r \in \Z_{\ge 0}$. 
As in the proof of Lemma \ref{Lemma 2.1}, the homogeneous part of degree 
\[
m_1\lambda_1+\cdots+ m_r\lambda_r+\mu s =
(m_1+\alpha_1 s)\lambda_1+\cdots+(m_r+\alpha_r s)\lambda_r
\]
in the equation (\ref{6}) yields
\begin{equation}
\xi^m\eta^s+c_1\xi^{m+\alpha}\eta^{s-1}+\cdots+c_s\xi^{m+s\alpha}=0, \label{7}
\end{equation}
where $c_1, \ldots,c_s \in k$ and 
$\xi^{m+i\alpha}=\xi_1^{m_1+i\alpha_1}\cdots\xi_r^{m_r+i\alpha_r}$ 
for $0 \le i \le s$. 
Hence $m_1=\cdots=m_r=0$ and the equation (\ref{7}) yields an equation 
\begin{equation}
\left(\frac{\eta}{\xi^\alpha}\right)^s
+c_1\left(\frac{\eta}{\xi^\alpha}\right)^{s-1}+\cdots+c_s=0. \label{8}
\end{equation}
Since $k$ is algebraically closed and 
the equation (\ref{8}) is a minimal equation, we have 
$\eta=c\xi_1^{\alpha_1}\cdots\xi_r^{\alpha_r}$ with $c \in k$. 
Hence $B=A=k[\xi_1,\cdots, \xi_r]$. Let $x_i=\xi_i$ for $1 \le i \le r$. 
Then $D$ is written as 
\[
D=\lambda_1 x_1\frac{\partial}{\partial x_1}+\cdots
+\lambda_r x_r\frac{\partial}{\partial x_r}.
\]
\QED

\smallskip

If $\lambda_1, \ldots, \lambda_r$ are not numerically independent, 
we are ignorant of any characterization theorem 
corresponding to Theorem \ref{Theorem 2.6}.

\section{Algebraic derivations and singularities}

Let $B$ be a normal affine domain of dimension two. If $D$ is a non-trivial locally nilpotent $k$-derivation  
on $B$, then $\Spec B$ has at most cyclic quotient singular points \cite{M2}. In the case of algebraic 
derivations, worse singularities can coexist as shown in the following result.

\begin{thm}\label{Theorem 3.1}
Let $G$ be a finite group acting on the $k$-vector space $kx_1+\cdots+kx_r$ linearly. Let $B$ be a 
polynomial ring $k[x_1,\ldots,x_r]$ and let $A=B^G$ the ring of $G$-invariants. Let 
$D=\sum_{i=1}^r x_i\frac{\partial}{\partial x_i}$ be the Euler derivation. Then $D$ induces a nontrivial algebraic 
derivation $D_A$ on $A$.
\end{thm}
\Proof
Note that $B$ is a graded ring with grading $\deg x_i=1$ for $1 \le i \le r$. Since $G$ acts linearly on $kx_1+ 
\cdots+kx_r$, $A$ is a graded subalgebra. Namely, $A=\bigoplus_{n\ge 0}A_n$, where $A_n=A\cap B_n$ with $B_n$ the 
set of homogeneous polynomials of degree $n$. It is then straightforward to see that $D(f)=nf$ for every $f \in A_n$.
Hence $D$ induces a $k$-derivation $D_A$ on $A$. Since $D$ is algebraic, so is $D_A$.
\QED

When $r=2$, the Euler derivation in Theorem \ref{Theorem 3.1} induces 
an algebraic derivation on $k[x,y]^G$, 
where $G$ is a finite subgroup of $\GL(2,k)$. 
However, we do not know the type of singularity 
which can coexist with algebraic derivations treated in Theorems 
\ref{Theorem 2.4} and \ref{Theorem 2.5}. 

\bigskip

\noindent
{\bf Acknowledgements}

The authors are grateful to Professor J. Winkelmann and Professor D. Daigle 
for pointing out that in Theorem \ref{Theorem 1.9} 
part of the condition on the totally orderedness of $M$ is redundant.

\end{document}